\documentstyle{amsppt}
\magnification=\magstep1
\NoRunningHeads
\NoBlackBoxes

\pageheight{9truein}
\pagewidth{6.5truein}
\TagsOnRight
\define\lra#1{{\left\langle#1\right\rangle}}
\define\xp#1{{\phantom{#1}}}
\define\({\left(}
\define\){\right)}
\define\[{\left[}
\define\]{\right]}

\let\en=\enspace

\define\ocal{{\Cal O}}
\define\codim{\operatorname{codim}}
\define\coker{\operatorname{coker}}
\define\grade{\operatorname{grade}}

\define\Syz{\operatorname{Syz}}

\define\Ext{\operatorname{Ext}}
\define\Hom{\operatorname{Hom}}

\redefine\Im{\operatorname{Im}}
\define\height{\operatorname{height}}

\define\rank{\operatorname{rank}}

\define\Tor{\operatorname{Tor}}

\topmatter
\title On a Consequence of the Order Ideal Conjecture \endtitle
\vglue-0.30truein
\author S.~P.\ Dutta \endauthor

\address\newline
University of Illinois\newline
Department of Mathematics\newline
1409 West Green Street\newline
Urbana, IL 61801\newline
U.S.A.\newline
\eightpoint E-mail: dutta\@math.uiuc.edu
\endaddress
%
\thanks
\hbox{\hskip-12pt}AMS Subject Classification: Primary 13D02, 13D22,
Secondary 13C15, 13D25, 13H05\newline
Key words and phrases: Projective dimension, grade, order ideal,
syzygy.\hfill\hfill 
\endthanks

\dedicatory Dedicated to H.-B.\ Foxby on his 65th birthday.\enddedicatory

\abstract
Given a minimal set of generators $\bold{x}$ of an ideal $I$ of height d in a
regular local ring ($R, m, k$) we prove several cases for which the map
$K_d(\bold{x}; R) \otimes k \to \Tor_d^R (R/I, k)$ is
the 0-map. As a consequence of the order ideal conjecture we derive several cases
for which $K_{d+i}(\bold{x}; R) \otimes k \to \Tor_{d+i}^R (R/I, k)$ are 0-maps
for $i \ge 0$.

\endabstract

\endtopmatter

\document

\vglue-18pt



\baselineskip18pt
\

In order to solve the syzygy problem in the equicharacteristic case
Evans and Griffith proved the following:

\proclaim{Theorem A {\rm (\cite{E-G1}; Th. 2.4, \cite{E-G2})}} Let $R$ be a
local ring containing a field. Let $M$ be a finitely generated $k$-th
syzygy of finite projective dimension and let $x$ be a minimal
generator of $M$. Then the order ideal $\ocal_M(x) = \{f(x)\mid f\in
\Hom_R(M,R)\}$ has grade at least $k$.
\endproclaim

We drop $M$ from the notation for order ideals when there is no scope
for confusion. Afterwards Bruns and Herzog extended the above theorem
for finite complexes of finitely generated free modules in the
following way.

\proclaim{Theorem B {\rm (Th.\ 9.5.2, \cite{B-H})}} Let $(R,m)$ be a
local ring containing a field and let
$$
F_\bullet : 0\to
F_s{\;\overset{\phi_s}\to{\longrightarrow}\;}F_{s-1}\to\cdots\to
F_1{\;\overset{\phi_1}\to{\longrightarrow}\;}F_0\to 0
$$
be a complex of finitely generated free $R$-modules. Then for
every $j$, $1\le j\le s$ and for every $e\in F_j$ with $e\not\in
mF_j+\Im\phi_{j+1}$, codimension $\ocal\(\phi_j(e)\)\ge \codim
F_\bullet + j$, where $\codim F_\bullet = \text{inf}\;
\{\text{codimension}\,I_{r_i}(\phi_i)-i\mid r_i=
\sum\limits^s_{j=i}(-i)^{j-i}\rank F_j\}$.
\endproclaim

The existence of big Cohen-Macaulay modules, due to Hochster ([H2]), played an
important role in the proofs of these two theorems.

Bruns and Herzog observed the following as a consequence of the above Theorems.

\proclaim{Result {\rm (Corollary 9.5.3, \cite{B-H})}} Let $(R,m,k)$ be
a regular local ring containing a field and $I\subset m$ be an ideal
generated by $x_1,\dots,x_s$. Let $K_\bullet ({\bold{x}};\,R)$ be the
corresponding Koszul complex. Then the natural map
$K_i({\bold{x}};\,R)\otimes k\to \Tor^R_i(R/I,\,k)$ is 0 for
$i>\grade I$.
\endproclaim

The question whether the statements of the Theorems and the Result mentioned above
are valid in mixed characteristic is very much open. A graded version of Theorem A
for mixed characteristic was proved later by Evans and Griffith
(\cite{E-G3}). These results induced us to propose the following conjecture (\cite{D1},
\cite{D2}).

\proclaim{Order ideal conjecture} Let $(R,m)$ be a local ring and let
$M$ be a finitely generated $k$-th syzygy of finite projective
dimension and let $x$ be a minimal generator of $M$. Then the order
ideal $\ocal_M(x)$ has grade at least~$k$.
\endproclaim

We say that a finitely generated R-module N of finite projective
dimension over R satisfies the order ideal conjecture if minimal
generators of all its syzygies do so.

In this note we would like to focus our attention to the observation made by Bruns
in the above Result. The important question that arises from our study of the
Monomial Conjecture,
due to Hochster ([H1]), for ideals $I$ of height $d$ in a regular local ring $R$ is the following: given
${\bold{x}}=\{x_1,\dots x_s\}$, a minimal set of generators of $I$, when is $K_d({\bold{x}};\,R)\otimes k\to
\Tor^R_d(R/I,\,k)$ the 0-map? In general the answer to the above
question is yet to be resolved for any characteristic. We have proved
in \cite{D1} that an affirmative answer to the above question for a
certain class of ideals implies the Monomial Conjecture. Thus, whether
$R$ is equicharacteristic or not, the vanishing of the above map is very significant.
We would mainly investigate this question in this note. Our main Result states the
following:

\proclaim{Theorem 1.2} Let $(R,m,k)$ be a regular local ring of dimension n
and let $I$ be an ideal of height d. We have the following:

\roster
\item"{a)}" Let ${\bold{x}}=\{x_1,\dots,x_s\}$ denote a minimal set of generators
of~$I$. If $R/I$ is $S_2$, then the natural map
$K_d({\bold{x}};\,R)\otimes k\to \Tor^R_d(R/I,\,k)$ is the 0-map.

\item"{b)}" Suppose that $R$ contains a field. Let ${\bold{\underline {x}}} \
\(resp. \ {\bold{x}}\)$ denote the ideal (resp. $R$-regular sequence ) generated by
$x_1,\dots,x_d$ in $I$. Let $S=R/{\bold{\underline{x}}}$ and let $\Omega$ denote the canonical module of $R/I$ i.e. $\Omega=\Hom_S \(R/I,\,S \) =\Ext^d_R(R/I,\,R)$.
If $\widetilde{\Omega}$ denotes the lift of $\Omega$ in $R$ via the natural surjection $R\to S$, then
$K_d({\bold{x}};\,R)\otimes k\to \Tor_d(R/\widetilde{\Omega},k)$ is
the $0$-map. Moreover, if $R/I$ is $S_2$, then for any set of minimal
generators $y_1,\dots,y_s$ of $\widetilde{\Omega}$,
$K_d(y_1,\dots,y_s;\,R)\otimes k\to \Tor_d(R/\widetilde{\Omega},\,k)$
is the $0$-map.

\item"{c)}" Let ${\bold{x}}$ denote a minimal set of generators of $I$. Then the natural maps
$K_{d+i}({\bold{x}};\,R)\otimes k\to\Tor_{d+i}(R/I,\,k)$ are
$0$-maps for $i>1$.

\item"{d)}" Assume that the order ideal conjecture is valid
for regular local rings of dimension (n-1). Let ${\bold{x}}$ denote a minimal set of generators of $I$.
If $I\cap(m - m^2) \ne \phi$ or depth
$R/I > 0 $ then the natural maps $K_{d+i}({\bold{x}};\,R)\otimes k\to \Tor^R_{d+i}({R/I},\,k)$
are $0$-maps for $i\ge 1$. Moreover if $R/I$ is $S_2$, then
$K_{d+i}({\bold{x}};\,R)\otimes k\to \Tor^R_{d+i}({R/I},\,k)$ are $0$-maps for $i\ge 0$

\item"" In particular, under the same hypothesis, for any prime ideal $P$ such that $R/P$ is
normal, $K_{d+i}({\bold{x}};\,R)\otimes k\to \Tor^R_{d+i}(R/P,\,k)$
are $0$-maps for $i\ge 0$.
\endroster
\endproclaim

\noindent We would like to mention that in \cite{D1} it was pointed out that the validity of
first half of part b) of the above theorem implies the Monomial Conjecture in mixed characteristic.

In the final observation (Prop. 1.3) of this note we prove the existence of
ideals $I$ of height $d$ in a regular local ring ($R, m, k$) for which
$K_i({\bold{x}}; R) \otimes k \to \Tor_i^R (R/I, k)$ is non-zero for $i \le d$;
here ${\bold{x}}$ denotes an R-sequence contained in $I$.

\head{Section 1}\endhead

We intend to prove our theorem in this section.

\proclaim{1.1\en Lemma}Let $(R,m)$ be a local ring such that the order ideal
conjecture is valid on $R$. Let $I$ be an ideal of $R$ of finite
projective dimension over $R$ and let $\grade I = d$. Let $x_1,\dots,
x_s$ be a set of generators of $I$ and let
$K_\bullet({\bold{x}};\,R)$ denote the corresponding Koszul
complex. Then the natural $K_i({\bold{x}},\,R)\otimes k\to
\Tor^R_i(R/I,\,k)$ is $0$ for $i\ge d+1$.
\endproclaim

\demo{Proof}The proof is a straightforward application of the validity
of the order ideal conjecture as in (9.5.3) \cite{B-H}. Let
$(F_\bullet,\,\phi_\bullet)$ be a minimal resolution of $R/I$ and let
$\psi_\bullet : K_\bullet ({\bold{x}};\,R)\to F_\bullet$ lift the
identity map in $R/I$. If for any minimal generator $e_{j_1}\wedge
e_{j_2}\wedge \cdots \wedge\,e_{j_i}$ of $K_i$, $\psi_i(e_{j_1}\wedge
\cdots \wedge \,e_{j_i})=e$ is a minimal generator of $F_i$, then by
commutativity  of the diagram constructed via $\psi_\bullet$,
$\ocal(\phi_i(e))$ has grade $\le d$. Due to our assumption on order
ideals, this can only happen if $i\le d$.
\enddemo

\demo{{\bf 1.2}\en Proof of Theorem \rm{\bf 1.2}}\hfill\newline

a)\en Let $z_1,\dots,z_d\in I$ be a regular sequence on $R$ and let
$S=R/(z_1,\dots,z_d)$. We will write ${\bold{\underline z}} (resp. \ {\bold {z}})$ to denote the
ideal (resp. sequence) generated by $z_1,\dots,z_d$. Let $\Omega=\Hom_S(R/I,\,S)(\simeq
\Ext^d_R(R/I,\,R))$ denote the canonical module of $R/I$. Let
$\widetilde{\Omega}$ denote the lift of $\Omega$ in $R$ via the
natural surjection $R\to S$. Since $R/I$ is $S_2$,
$\Hom_S(\Omega,\,\Omega)\simeq R/I$. Hence, by applying
$\Hom_S(-,\,S)$ to the short exact sequence
$$
0\to \Omega\to S\to S/\Omega\to 0
$$
we obtain $\Hom_S(S/\Omega,\,S)=I/{\bold{\underline z}}$ and
$\Ext^1_S(S/\Omega,\,S)=0$. Since $S/\Omega\simeq
R/\widetilde{\Omega}$ and height of $I$ is d, we have
$$
\Ext^d_R(R/\widetilde{\Omega},\,R) \simeq \Hom_S (R/\widetilde{\Omega}, S) =
\Hom_S (S/\Omega, S)=I/{\bold{\underline z}}
$$
and
$$
\Ext^{d+1}_R(R/\widetilde{\Omega},\,R)= \Ext_S^1 (R/\widetilde {\Omega}, S) = \Ext_S^1 (S/\Omega, S) = 0.\tag{1}
$$

Let $F_\bullet,\,P_\bullet$ denote the minimal free resolutions of
$R/I$ and $I/({\bold{\underline z}})$ respectively and let
$K_\bullet^\prime=K_\bullet(z_1,\dots, z_d;\,R)$ be the Koszul complex
corresponding to the sequence $z_1,\dots, z_d$. Let $\psi_\bullet :
P_\bullet \to K_\bullet^\prime$ denote a lift of the natural injection
$I/{\bold{\underline z}}\to S=R/{\bold{\underline z}}$ and let $\phi_\bullet : K_\bullet^\prime \to
F_\bullet$ denote a lift of the natural surjection $S=R/{\bold{\underline z}}\to
R/I$.
\enddemo

\proclaim{Claim}i)\en $\psi_d(P_d)=K_d^\prime=R$ and ii)\en
$\phi_d(K_d^\prime)\subset m\,F_d$.
\endproclaim

\demo{Proof of Claim i)} Let $(L_\bullet,\,\alpha_\bullet)$ be a
minimal free resolution of $R/\widetilde{\Omega}$ over $R$ and
$\theta_\bullet : K_\bullet^\prime \to L_\bullet$ be a lift of the natural
surjection $R/{\bold{\underline z}}\to R/\widetilde{\Omega}$. Let $L^\ast_\bullet
= \Hom_R(L_\bullet,\,R)$. We consider the following commutative
diagram
$$
{\eightpoint \matrix
0 & \to & R & \to & L_1^\ast & \to & L_2^\ast & \to & \cdots & \to &
L_d^\ast  \to & L^\ast_{d+1}  \to & L^\ast_{d+2}  \to & G_{d+2} & \to & 0\\
& & \| & & \theta^\ast_1\Big\downarrow\quad & &
\theta^\ast_2\Big\downarrow\quad & & & & \theta^\ast_d\Big\downarrow\quad\\
0 & \to & R & \to & R^d & \to & R^{d_2} & \to & \cdots & \to & R
\to & R/{\bold{\underline z}}  \to  0
\endmatrix
}\tag{2}
$$
Here $G_{d+i}=\coker \alpha^\ast_{d+i}$, for $0\le i\le 2$. Since
height $\widetilde{\Omega}$ is $d$, we have $H^i(L^\ast_\bullet)=0$ for $i<d$ and
$H^d(L^\ast_\bullet)=\Ext^d(R/\widetilde{\Omega},\,R)=I/{\bold{\underline z}}$
$\hookrightarrow G_d$. Also it follows from (1) that
$H^{d+1}(L^\ast_\bullet) = \Ext^{d+1}_R(R/\widetilde{\Omega},\,R)=0$.
Hence $G_{d+1}$ is a submodule of $L^\ast_{d+2}$ and
$$
L^\prime_\bullet : 0\to R\to L^\ast_1\to \cdots\to L^\ast_d\to 0
$$
is a minimal free resolution of $G_d$. Moreover, the map $G_d \to R/{\bold{\underline z}}$,
induced by $\theta^\ast_d$ in (2),
induces the natural inclusion $I/{\bold{\underline z}}\hookrightarrow
R/{\bold{\underline z}}$. Let $\beta_\bullet : P_\bullet \to L^\prime_\bullet$ be
a lift of the natural inclusion $I/{\bold{\underline z}}\hookrightarrow G_d$.
We have the following commutative diagram:
$$
{\eightpoint
\matrix & \to & P_d & \to & P_{d-1} & \to & \cdots & \to & P_0 & \to &
I/({\bold{\underline z}}) & \to & 0\\
& & \quad\Big\downarrow \beta_d& & \Big\downarrow & & & & \quad\Big\downarrow\beta_0 & &
\Big\downarrow\\
0 & \to & R & \to & L^\ast_1\en & \to & \cdots & \to & L_d^\ast\en & \to & G_d
\endmatrix
}\tag{3}
$$
\enddemo

\proclaim{Subclaim}$\beta_d(P_d)=R$.
\endproclaim

\demo{Proof of the Subclaim}If possible, let $\beta_d(P_d)\subset
m$. Since $\coker(I/{\bold{\underline{z}}} \to G_d)=\Im \alpha^\ast_{d+1}$, the
mapping cone of $\beta_\bullet$ gives rise to a free resolution of
$G_{d+1}$. Let $(Q_\bullet,\,\gamma_\bullet)$ be a minimal free
resolution of $G_{d+1}$ extracted from this mapping cone. Since
$\beta_d(P_d)\subset m$, the copy of $R=L_0{^\prime}$ is a summand of $Q_{d+1}$;
we denote this copy by $Re$. Since height $\widetilde{\Omega}$ is $d$ and $L_\bullet$ is
a minimal free resolution of $R/\widetilde\Omega$ , we have

$$
  {\height}\  {\ocal}  (\gamma_{d+1} (e))\le d \ \ (diagram (3)). \tag{4}
$$

If $R$ is equicharacteristic, this cannot happen due to Theorem A. Let
us assume that $R$ is of mixed characteristic $p>0$. Since $G_{d+1}$ is a
submodule of $L^\ast_{d+2}$, $p$ is a non-zero-divisor on
$G_{d+1}$. Let $\overline{G}_{d+1}=G_{d+1}\otimes R/pR$ and
$(\overline{Q}_\bullet,\,\overline{\gamma}_\bullet)=(Q_\bullet\otimes
R/pR,\,\gamma_\bullet\otimes R/pR)$. Then
$(\overline{Q}_\bullet,\,\overline{\gamma}_\bullet)$ is a minimal free
resolution of $\overline{G}_{d+1}$ over $\overline{R}(=R/pR)$. Since
$\overline{R}$ is an equicharacteristic local ring (characteristic
$p>0$), by Theorem A, height $\ocal(\gamma_{d+1}(\overline{e}))$ must
be $\ge d+1$. This contradicts (4) and hence $\beta_d(P_d)=R$.
\enddemo

Since $\psi_\bullet$ can be constructed as the composite of $\theta_\bullet^*$
and $\beta_\bullet$, the above observation implies claim
i)~i.e.\ $\psi_d(P_d)=K_d^\prime=R$ via combination of diagrams~(2) and (3).

\demo{Proof of claim ii)}Since $\psi_d(P_d)=R$, in the minimal free
resolution $F_\bullet$ of $R/I$ extracted from the mapping cone of
$\psi_\bullet$, the copy of $K_d^\prime=R$ does not survive as a minimal
generator of $F_d$ i.e.\ $\phi_d(K_d^\prime)\subset m\,F_d,$ where
$\phi_\bullet$ is extracted from the natural inclusion of $K_\bullet^\prime$
into the mapping cone of $\psi$.
\enddemo

Thus both parts of our claim are proved. Since height $I$ is $d$ and $R$ is
regular local, grade $I$ is also $d$. We can find a minimal set of generators
$x_1,\dots,x_s$ of $I$ such that any subset $x_{i_1},\dots,x_{i_d}$,
$1\le i_1<i_2<\cdots<i_d\le s$, form an $R$-sequence. Let
$K_\bullet=K_\bullet(x_1,\dots,x_s;\,R)$ denote the Koszul complex
corresponding to $x_1,\dots,x_s$. Let $F_\bullet$ be a minimal free
resolution of $R/I$ over $R$ and let $\phi_\bullet : K_\bullet\to
F_\bullet$ denote a lift of the identity map on $R/I$. By part ii) of
our claim, $\phi_d(R\,e_{i_1},\wedge\cdots \wedge\,e_{i_d})\subset
m\,F_d,\, \forall d$-tuple $(i_1,\dots,i_d)$ such that $1\le i_1<i_2<\cdots
<i_d\le s$. Hence $K_d({\bold{x}};\,R)\otimes R/m\to F_d\otimes R/m$
is the $0$-map.

b)\en Let $S = R/{\bold{\underline x}}$ and let
$\Omega=\Hom_S(R/I,\,S)= \Ext^d_R(R/I,\,R)$. Let
$(P_\bullet,\,\alpha_\bullet)$ be a minimal free resolution of $R/I$
over $R$; let $P^\ast_\bullet = \Hom_R(P,\,R)$. We consider the
following complex truncated from $P^\ast_\bullet$:
$$
0\to R\to P_1^\ast\to \cdots\to P_d^\ast\to P^\ast_{d+1}\to G_{d+1}\to
0\tag{1}
$$
where $G_{d+1}=\coker \alpha^\ast_{d+1}$ and $G_d=\coker
\alpha^\ast_d$. We have $H^i(P^\ast_\bullet)=0$ for $i<d$ and
$H^d(P^\ast_\bullet) = \Omega$. Let $P^\prime_\bullet:0\to R\to
P^\ast_1\to \cdots \to P_d^\ast \to 0$ denote the minimal free
resolution of $G_d$ obtained from (1). Let $F_\bullet$ be a minimal
free resolution of $\Omega$ over $R$ and let $\phi_\bullet :
F_\bullet\to P^\prime_\bullet$ lift the natural inclusion
$\Omega\hookrightarrow G_d$. $R$ being equicharacteristic, Theorem A is
valid in $R$. Since the mapping cone of $\phi_\bullet$ gives rise to a
free resolution of $G_{d+1}$, arguing as in the proof of part a), it
follows that $\phi_d(F_d)=R$. Let $\theta_\bullet$:$K_\bullet ({\bold{x}};\,R) \to P_\bullet$
denote a lift of the natural surjection $R/\underline{\bold{x}} \to R/I$
and let $\theta_\bullet^* : P_\bullet^{\prime} \to K_\bullet ({\bold{x}};\,R)$
denote the dual map that lifts the map $G_d \to R/{\underline{\bold{x}}}$ induced by
$\theta_d^*$. Then $\theta_\bullet^*.\phi_\bullet : F_\bullet \to K_\bullet({\bold{x}}; R)$
lifts the natural injection $\Omega \hookrightarrow R/\underline{\bold{x}}$
(composition of $\Omega \hookrightarrow G_d \to R/\underline {\bold{x}}$) and
$\theta_0^* \phi_d (F_d) = R$. Let $L_\bullet$ be a minimal free
resolution of $R/\widetilde{\Omega}$ extracted from the mapping cone of
$\theta_\bullet ^*. \phi_\bullet$ and $\psi_\bullet:
K_\bullet({\bold{x}};\,R)\to L_\bullet$ be a lift of the natural
surjection $R/\underline{\bold{x}}\to R/\widetilde{\Omega}$ induced by the inclusion map
of $K_\bullet({\bold x}; R)$ into the mapping cone of $\theta_\bullet^*. \phi_\bullet$. Again arguing as in
the proof of part a) we obtain $\psi_d(K_d ({\bold{x}}; R) = R) \subset m\,L_d$. Hence
the natural map $K_d({\bold{x}};\,R)\otimes k\to \Tor_d^R(R/\widetilde{\Omega},\,k)$
is the $0$-map.

Now assume $R/I$ is $S_2$. Then
$\Ext_R^{d+1}(R/\widetilde{\Omega},\,R)=0$ ((1) in part a)). Let ${\bold{y}} =  \{y_1,\dots,y_d\}$ be a
part of a minimal set of generators of $\widetilde{\Omega}$ such that
they form an $R$-sequence. Let $T=R/(y_1,\dots,y_d)$, and
$\widetilde{\Omega}/(y_1,\dots,y_d) = \omega$; then $T/\omega =
R/\widetilde{\Omega}$ and $\Ext^1_T(T/\omega,\,T) = \Ext_R^{d+1} (R/\widetilde\Omega, R) = 0$.

Let $J$ = $\Hom_{T}(T/\omega,\,T)$ = $\Ext^d_R(R/\widetilde{\Omega},\,R)$. We
consider the short exact sequence
$$
0 \to \omega\to T\to T/\omega\to 0.
$$
Applying $\Hom_T(\underline{\enspace},\,T)$ to this sequence, we obtain
$$
0 \to J\to T\to \Hom_T(\omega,\,T)\to 0.
$$
Hence $T/J\simeq \omega^\ast = \Hom_T(\omega,\,T)$. Since $T/\omega = R/\widetilde\Omega$
is $S_1$, $\omega$ is $S_2$ and hence $\Hom_T(T/J,\,T)=\omega$.

Let $\widetilde{J}$ be a lift of $J$ in $R$ via the natural surjection
$R\to T$. Let $K_\bullet({\bold {y}}; \ R)$ denote the
Koszul complex corresponding to $y_1,\dots,y_d$. Let $B_\bullet$ be a
minimal free $R$-resolution of $\omega$ and let $L_\bullet$ be a
minimal free resolution of $T/\omega=R/\widetilde{\Omega}$ over
$R$. Then, arguing as in part a), starting with a minimal free
resolution of $R/\widetilde{J}$ over $R$ we can show that the
natural map $K_d({\bold{y}};\,R)\otimes k\to
\Tor_d^R(R/\widetilde{\Omega},\,k)$ is the $0$-map. Now considering a minimal
set of generators $y_1,\dots,y_s$ of $\widetilde{\Omega}$ such for any
$d$-tuples $(i_1,\dots,i_d)$, $1\le i_1<i_2<\cdots<i_d\le s$,
$y_{i_1},\dots,y_{i_d}$ from an $R$-sequence, it follows that the
natural map $K_d(y_1,\dots,y_s;\,R)\otimes k\to
\Tor_d^R(R/\widetilde{\Omega},\,k)$ is the $0$-map.

c)\en This follows from the facts that for a minimal free
resolution $F_\bullet$ of a finitely generated module $M$ over $R$, the $(d+1)$th
syzygy $\Syz_R^{d+1}(M)=\Syz_R^d(S_1)$ where $S_1=\Syz_R^1(M)$ and the
mixed characteristic $p$ is a non-zero-divisor on $S_1$. Let
$\overline{R}=R/pR$. Then $\overline{F}_\bullet = F_\bullet(+1)\otimes
\overline{R}$ is a minimal free resolution of
$\overline{S}_1=S_1\otimes\overline{R}$ over $\overline{R}$. Since
$\overline{R}$ is equicharacteristic, Theorem A is valid on
$\overline{R}$ and hence our assertion follows.

d)\en If w $\epsilon (m-m^2)$ is a non-zero-divisor on $R/I$, the result follows by
similar arguments as in the proof of part c) above.

If w $\epsilon (m-m^2) \cap I$, the result follows from the proof of Lemma (1.1)
and the following theorem which is a part of theorem 2.5 in~\cite{D2}.

\proclaim{Theorem {\rm (Theorem 2.5, \cite{D2})}}Let $(R, m)$ be a regular local ring
and let $M$ be a module of finite projective dimension over $R$. Assume that the
order ideal conjecture is valid for regular local rings of dimension (n-1).
If $ann_RM\cap(m-m^2) \ne \phi$ or $depth_RM>0$, then M satisfies the order ideal conjecture.

\endproclaim
For a proof we refer the reader to theorem 2.5 in ([D2]).

The last part of the statement in d) now follows from a) and the above theorem.

{\bf 1.3.} In the final proposition we prove the existence of almost
complete intersection ideals I of height d for which $K_i({\bold{x}};\,R)\otimes k\to
\Tor^R_i(R/I,\,k)$ is non-zero for $0 \le i \le d $.

\proclaim {Proposition} Let (R, m, k) be a regular local ring and let P be a prime ideal of height d. Then there exists a family of almost complete intersection ideals $\{I_t\}_{t\epsilon N} \subset P$ such that $K_i({\bold {x}};\,R)\otimes k\to
\Tor^R_i(R/I_t,\,k)$ is injective for $i\ge0$ and $t>1$; here ${\bold{x}}$ denotes a part of a minimal
set of generator $x_1, \dots, x_d$ of $I_t$ for all t.
\endproclaim

\demo {Proof} We choose $x_1, \dots, x_d \in P$ such that $x_1, \dots, x_d$ form an
R-sequence and $PR_P = {\bold{\underline{x}}}R_P$; ${\bold{\underline x}}$ = the ideal generator by $x_1, \dots, x_d$.
We have a primary decomposition of ${\bold{\underline x}}$ as follows:

$$
{\bold {\underline x}} = P \cap q_2 \cap \dots \cap q_r
$$
where $q_i$ is $P_i$-primary, $2 \le i \le r$. Let $\Omega$ denote the canonical module for $R/P$;
then $\Omega = \Ext^d_R (R/P, R)$ = $\Hom_R(R/P, R/{\bold{\underline x}}) =
q_2 \cap \dots \cap q_r/{\bold{\underline x}}$. Let $\lambda \in P- \cup P_i$. It can be easily checked
that $\Omega = \Hom_R(R/({\bold{\underline x}}, \lambda^t), R/{\bold{\underline x}})$ for $t>0$.
Let $\Omega^\prime$ denote the lift of $\Omega$ in R
via the natural surjection $R \to R/{\bold{\underline x}}$.
Then $({\bold{\underline x}} + \lambda^tR)/\lambda^tR \simeq {\bold{\underline x}}/
{\bold{\underline x}} \cap \lambda^tR = {\bold{\underline x}}/\lambda^t \Omega^\prime$ for $t>0$.
Consider the injection $i_t : \Omega^\prime \to {\bold{\underline x}}$ defined by
$i_t(\mu) = \lambda^t\mu$. Let $F_\bullet$ denote a minimal free resolution of
$\Omega^\prime$, $K_\bullet = K_\bullet({\bold {x}}; R) (+1)$ (i.e. $K_i = K_{i+1}
({\bold {x}}; R)$) denote the minimal free resolution of ${\bold{\underline x}}$
and let $\Phi_\bullet$ : $F_\bullet \to K_\bullet$ denote a lift of $i_t$. Since
$\lambda \centerdot i_t = i_{t+1}, \lambda \centerdot \Phi_\bullet$ is a lift of $i_{t+1}$. Let us assume
that $t>1$. Then im$\Phi_j \subset mK_j$ for $j\ge0$. Hence $L_\bullet$, the mapping cone of
$\Phi_\bullet$, is a minimal free resolution of ${\bold{\underline x}}/\lambda^t\Omega^\prime$.
Hence, via the natural injection of $K_\bullet$ into $L_\bullet$,
$K_\bullet$ imbeds into $L_\bullet$ i.e. $K_j$ is a free summand of $L_j$ for $j \le d-1$.
Let $K_\bullet(\lambda^t; R)$ be the Koszul complex corresponding to $\lambda^t$ over R.
Let $\Psi_\bullet : L_\bullet \to K_\bullet (\lambda^t; R)$ lift the inclusion map :
${\bold{\underline{x}}}/\lambda^t\Omega^\prime \hookrightarrow R/\lambda^tR$. Then the mapping
cone $P_\bullet$ of $\Psi_\bullet$ is a minimal free resolution of
$R/({\bold{\underline {x}}},\lambda^t)$ and it is easy to check that
$K_\bullet({\bold {x}}; R)$ imbeds into $P_\bullet$ i.e.
$K_i({\bold{\underline {x}}}; R)$ is a free summand of $P_i$ for $i \ge 0$. Hence
$K_i({\bold{\underline {x}}}; R) \otimes k \to \Tor^R_i(R/({\bold{\underline {x}}},
\lambda^t), k)$ is injective for $i\ge0$ and our proof is complete.
\enddemo


\enddocument